\def\TheMagstep{\magstep1}	% Normal magnification
		% Changed to \magstep0 by \DoublepageOutput{TRUE} 
\def\PaperSize{letter}		% \PaperSize is used to

\magnification=\magstep1

\let\:=\colon  
   \let\?=\overline

\let\Sum=\sum \def\sum{\Sum\nolimits}
\def\Bl{\rm{ Bl}}
\def\IC{{\bf C}} 
\def\IP{{\bf P}}

\def\Pd #1#2{{\partial#1\over\partial#2}}
	
\def\and{\hbox{ and }}
	\def\Af{\hbox{\rm A$_f$}}

\def\DONE{*!*}
\def\NextDef #1 {\def\NextOne{#1}%
 \ifx\NextOne\DONE\let\next\relax
 \else\expandafter\xdef\csname#1\endcsname{\TheOp}
  \let\next\NextDef
 \fi \next}
\def\TheOp{\mathop{\rm\NextOne}}
 \NextDef 
  Projan Supp Proj Sym Spec Hom cod Ker dist
 *!*
\def\TheOp{{\cal\NextOne}}
\NextDef 
  E F G H I J M N O R S
 *!*
\def\TheOp{\hbox{\rm\NextOne}}
\NextDef 
 A ICIS 
 *!*
 
 %%  STYLE MACROS  %%
%% Redefine \item to give greater indentation than AMSTeX
%   and the roman font within parentheses.   
\def\item#1 {\par\indent\indent\indent\indent \hangindent4\parindent
 \llap{\rm (#1)\enspace}\ignorespaces}
%% Define a similar macro without the hanging indentation for assertions
%% and that starts each part with an ordinary \parindent
 \def\inpart#1 {{\rm (#1)\enspace}\ignorespaces}
 \def\part {\par\inpart}
%% For displaying statements in italics with narrower margins and numbers
%\def\state{\smallskip\begingroup\narrower\noindent\it}

\catcode`\@=11		% make @ a letter temporarily

%% Modification of the PLAIN footnote macro for 8pt
\def\vfootnote#1{\insert\footins\bgroup
 \eightpoint %% only change
 \interlinepenalty\interfootnotelinepenalty
  \splittopskip\ht\strutbox % top baseline for broken footnotes
  \splitmaxdepth\dp\strutbox \floatingpenalty\@MM
  \leftskip\z@skip \rightskip\z@skip \spaceskip\z@skip \xspaceskip\z@skip
  \textindent{#1}\footstrut\futurelet\next\fo@t}

%%  ``Ties'' with a \thinspace for page numbers
\def\p.{p.\penalty\@M \thinspace}
\def\pp.{pp.\penalty\@M \thinspace}
%% For Roman parenthetical material in nonRoman text
\def\(#1){{\rm(#1)}}\let\leftp=(
\def\activeleftp{\catcode`\(=\active}
{\activeleftp\gdef({\ifmmode\let\next=\leftp \else\let\next=\(\fi\next}}

%% Sectioning
\def\sct#1\par
  {\removelastskip\vskip0pt plus2\normalbaselineskip \penalty-250 
  \vskip0pt plus-2\normalbaselineskip \bigskip
  \centerline{\smc #1}\medskip}

\newcount\sctno \sctno=0
\def\sctn{\advance\sctno by 1 
% {\bf\hbox to \parindent{\number\sctno.\hfil}#1.}\enspace\ignorespaces}
 \sct\number\sctno.\quad\ignorespaces}

%% Display numbers
\def\dno#1${\eqno\hbox{\rm(\number\sctno.#1)}$}
\def\Cs#1){\unskip~{\rm(\number\sctno.#1)}}

%% For setting results
\def\proclaim#1 #2 {\medbreak
  {\bf#1 (\number\sctno.#2)}\enspace\bgroup\activeleftp
\it}
\def\endproclaim{\par\egroup\medskip}
\def\pf{\endproclaim{\bf Proof.}\enspace}
\def\lem{\proclaim Lemma } \def\prp{\proclaim Proposition }
\def\cor{\proclaim Corollary }	\def\thm{\proclaim Theorem }
\def\rmk#1 {\medbreak {\bf Remark (\number\sctno.#1)}\enspace}
\def\eg#1 {\medbreak {\bf Example (\number\sctno.#1)}\enspace}

%% PAGE LAYOUT
\parskip=0pt plus 1.75pt \parindent10pt
\hsize29pc
\vsize44pc
\abovedisplayskip6pt plus6pt minus2pt
\belowdisplayskip6pt plus6pt minus3pt

\def\TRUE{TRUE}	% For Boolean tests
\ifx\DoublepageOutput\TRUE \def\TheMagstep{\magstep0} \fi
\mag=\TheMagstep

% CENTER TEXT ON PAGE
	% additional vertical adjustment
\newskip\vadjustskip \vadjustskip=0.5\normalbaselineskip
\def\centertext
 {\hoffset=\pgwidth \advance\hoffset-\hsize
  \advance\hoffset-2truein \divide\hoffset by 2\relax
  \voffset=\pgheight \advance\voffset-\vsize
  \advance\voffset-2truein \divide\voffset by 2\relax
  \advance\voffset\vadjustskip
 }
\newdimen\pgwidth\newdimen\pgheight
\def\letter{letter}\def\AFour{AFour}
\ifx\PaperSize\letter
 \pgwidth=8.5truein \pgheight=11truein 
 \message{- Got a paper size of letter.  }\centertext 
\fi
\ifx\PaperSize\AFour
 \pgwidth=210truemm \pgheight=297truemm 
 \message{- Got a paper size of AFour.  }\centertext
\fi

%% TWO-COLUMN LANDSCAPE FORMAT
% Modified from the TeX book, p. 257.
 \newdimen\fullhsize \newbox\leftcolumn
 \def\fulline{\hbox to \fullhsize}
\def\doublepageoutput
{\let\lr=L
 \output={\if L\lr
          \global\setbox\leftcolumn=\columnbox \global\let\lr=R%
        \else \doubleformat \global\let\lr=L\fi
        \ifnum\outputpenalty>-20000 \else\dosupereject\fi}%
 \def\doubleformat{\shipout\vbox{%
        \fulline{\hfil\hfil\box\leftcolumn\hfil\columnbox\hfil\hfil}%
				}%
		  }%
 \def\columnbox{\vbox
   {\makeheadline\pagebody\makefootline\advancepageno}%
   }%
 \fullhsize=\pgheight \hoffset=-1truein
 \voffset=\pgwidth \advance\voffset-\vsize
  \advance\voffset-2truein \divide\voffset by 2
  \advance\voffset\vadjustskip
 \let\firstheadline=\hfil
 
% \null\vfill\nopagenumbers\eject\pageno=1\relax % to put page on right
}
\ifx\DoublepageOutput\TRUE \doublepageoutput \fi

%% ADDITIONAL FONTS
 \font\twelvebf=cmbx12		% For title
 \font\smc=cmcsc10		% For authors' names

%% EIGHT POINT TYPE FOR FOOTNOTES AND REFERENCES
\def\eightpoint{\eightpointfonts
 \setbox\strutbox\hbox{\vrule height7\p@ depth2\p@ width\z@}%
 \eightpointparameters\eightpointfamilies
 \normalbaselines\rm
 }
\def\eightpointparameters{%
 \normalbaselineskip9\p@
 \abovedisplayskip9\p@ plus2.4\p@ minus6.2\p@
 \belowdisplayskip9\p@ plus2.4\p@ minus6.2\p@
 \abovedisplayshortskip\z@ plus2.4\p@
 \belowdisplayshortskip5.6\p@ plus2.4\p@ minus3.2\p@
 }
\newfam\smcfam
\def\eightpointfonts{%
 \font\eightrm=cmr8 \font\sixrm=cmr6
 \font\eightbf=cmbx8 \font\sixbf=cmbx6
 \font\eightit=cmti8 
 \font\eightsmc=cmcsc8
 \font\eighti=cmmi8 \font\sixi=cmmi6
 \font\eightsy=cmsy8 \font\sixsy=cmsy6
 \font\eightsl=cmsl8 \font\eighttt=cmtt8}
\def\eightpointfamilies{%
 \textfont\z@\eightrm \scriptfont\z@\sixrm  \scriptscriptfont\z@\fiverm
 \textfont\@ne\eighti \scriptfont\@ne\sixi  \scriptscriptfont\@ne\fivei
 \textfont\tw@\eightsy \scriptfont\tw@\sixsy \scriptscriptfont\tw@\fivesy
 \textfont\thr@@\tenex \scriptfont\thr@@\tenex\scriptscriptfont\thr@@\tenex
 \textfont\itfam\eightit	\def\it{\fam\itfam\eightit}%
 \textfont\slfam\eightsl	\def\sl{\fam\slfam\eightsl}%
 \textfont\ttfam\eighttt	\def\tt{\fam\ttfam\eighttt}%
 \textfont\smcfam\eightsmc	\def\smc{\fam\smcfam\eightsmc}%
 \textfont\bffam\eightbf \scriptfont\bffam\sixbf
   \scriptscriptfont\bffam\fivebf	\def\bf{\fam\bffam\eightbf}%
 \def\rm{\fam0\eightrm}%
% \tt \ttglue=0.5em plus0.25em minus0.15em
 }

%% HEADLINE STYLE
\def\today{\ifcase\month\or	% From the TeX book p. 406
 January\or February\or March\or April\or May\or June\or
 July\or August\or September\or October\or November\or December\fi
 \space\number\day, \number\year}
\nopagenumbers
\headline={%
  \ifnum\pageno=1\firstheadline
  \else
    \ifodd\pageno\oddheadline
    \else\evenheadline\fi
  \fi
}
\let\firstheadline\hfill
\def\oddheadline{\eightpoint \rlap{\today}
 \hfil\headtitle\hfil\llap{\folio}}
\def\evenheadline{\eightpoint\rlap{\folio}
 \hfil\author\hfil\llap{\today}}
\def\headtitle{\title}

%% REFERENCING
	% to introduce the keys in order
 \newcount\refno \refno=0	 \def\NoKey{*!*}
 \def\MakeKey{\advance\refno by 1 \expandafter\xdef
  \csname\TheKey\endcsname{{\number\refno}}\NextKey}
 \def\NextKey#1 {\def\TheKey{#1}\ifx\TheKey\NoKey\let\next\relax
  \else\let\next\MakeKey \fi \next}
 \def\RefKeys #1\endRefKeys{\expandafter\NextKey #1 *!* }
	% to set references
\def\SetRef#1 #2,#3\par{%
 \hang\llap{[\csname#1\endcsname]\enspace}%
  \ignorespaces{\smc #2,}
  \ignorespaces#3\unskip.\endgraf
 }
 \newbox\keybox \setbox\keybox=\hbox{[8]\enspace}
 \newdimen\keyindent \keyindent=\wd\keybox
\def\references{%\vskip-\smallskipamount
  \bgroup   \frenchspacing   \eightpoint
   \parindent=\keyindent  \parskip=\smallskipamount
   \everypar={\SetRef}}
\def\endreferences{\egroup}

%% SERIALS
 \def\serial#1#2{\expandafter\def\csname#1\endcsname ##1 ##2 ##3
  {\unskip\ #2 {\bf##1} (##2), ##3}}
 \serial{ajm}{Amer. J. Math.}
  \serial {aif} {Ann. Inst. Fourier}
 \serial{asens}{Ann. Scient. \'Ec. Norm. Sup.}
 \serial{comp}{Compositio Math.}
 \serial{conm}{Contemp. Math.}
 \serial{crasp}{C. R. Acad. Sci. Paris}
 \serial{dlnpam}{Dekker Lecture Notes in Pure and Applied Math.}
 \serial{faa}{Funct. Anal. Appl.}
 \serial{invent}{Invent. Math.}
 \serial{ma}{Math. Ann.}
 \serial{mpcps}{Math. Proc. Camb. Phil. Soc.}
 \serial{ja}{J. Algebra}
 \serial{splm}{Springer Lecture Notes in Math.}
 \serial{tams}{Trans. Amer. Math. Soc.}

	% modified \cite code from AMSTeX
\def\UThin{\penalty\@M \thinspace\ignorespaces}
	% unbreakable \thinspace for use after periods
\def\relaxnext@{\let\next\relax}
\def\cite#1{\relaxnext@
 \def\nextiii@##1,##2\end@{\unskip\space{\rm[\SetKey{##1},\let~=\UThin##2]}}%
 \in@,{#1}\ifin@\def\next{\nextiii@#1\end@}\else
 \def\next{{\rm[\SetKey{#1}]}}\fi\next}
\newif\ifin@
\def\in@#1#2{\def\in@@##1#1##2##3\in@@
 {\ifx\in@##2\in@false\else\in@true\fi}%
 \in@@#2#1\in@\in@@}
\def\SetKey#1{{\bf\csname#1\endcsname}}

\catcode`\@=12 %\active  %at signs are no longer letters

\def\title{ Invariants of $D(q,p)$ singularities}
\def\author{Terence Gaffney}
\RefKeys  B B-G F  G-2  G-3 G-G  GM   
LJT  L-T Mac  Ma Ma1  Ma2 P P1  S1 S2  
 \endRefKeys

\def\topstuff{\leavevmode
 \bigskip\bigskip
 \centerline{\twelvebf \title}
 \bigskip
 \centerline{\author}
 \medskip\centerline{\today}
\bigskip\bigskip}
\topstuff
\sct Introduction

The study of the geometry of non-isolated hypersurface singularities was begun by Siersma and his students (\cite{S1},\cite{S2},\cite{P},\cite{P1}).
 The basic examples of such functions defining these singularities 
are the 
$A(d)$ singularities and the $D(q,p)$ singularities. The $A(d)$ singularities, up to analytic equivalence, are the product of a Morse function and 
the zero map, 
while the simplest $D(q,p)$ singularity is the Whitney umbrella. These are the basic examples, because they correspond to stable germs of functions in 
the study of germs of functions with 
isolated singularities.  Given a germ of a function which defines a  non-isolated hypersurface singularity at the origin, which in the appropriate sense,
 has finite codimension in the set of such germs, the singularity type
of such germs away from the origin is $A(d)$ or $D(q,p)$. However, some of the basic invariants of the germs of type $D(q,p)$ have not been calculated yet. 
In this note we calculate the homotopy type of the Milnor fiber 
of germs of type $D(q,p)$, as well as their L\^e numbers. The calculation of the L\^e numbers involves the 
use of an incidence variety 
which may be useful for studying germs of finite codimension. The calculation shows that the set of symmetric matrices of kernel rank $\ge 1$.
  is an example of a hypersurface singularity with 
a Whitney stratification (given by the rank of the matrices) in which only one singular stratum gives 
a component of top dimension of the singular set of the conormal.

The results of this paper can be applied to the study of any function in which the generic singularity type is $A(d)$ or $D(p,q)$, thus to all the germs which are finitely determined in the sense of Pellikaan. This is because in the case of non-isolated singularities, the geometry of many strata may contribute to invariants at the origin. We describe such an application, which appeared in the recent paper of de Bobadilla and Gaffney (\cite{B-G}), which computes the Euler invariant at the origin of the zero set of a function which has only singularities of type $A(k)$ or $D(q,p)$ off the origin.

\sctn The Milnor fiber and the L\^e numbers of the $D(q,p)$ singularities

If $f$ has a non-isolated singularity of type $D(p(p+1)/2,p)$, then it is known that by a change of coordinates, $f$ has the normal form

$$f({\bf x},{\bf y})=(\Sum_{i\le j}^p x_{i,j}y_iy_j)+y^2_{p+1}+\dots+y^2_{p+k}=[{\bf y}]^t[{\bf X}][{\bf y}]+y^2_{p+1}+\dots+y^2_{p+k}.$$

Here $[{\bf X}]$ is a symmetric matrix with diagonal entries $x_{i,i}$ and off diagonal entries $1/2x_{i,j}$ and $n$, the dimension of the domain of $f$ is $p+k+p(p+1)/2$. 
 In general, in the notation $D(q,p)$, $p$ refers to the size of the matrix $[{\bf X}]$ while the number of generators of the ideal $I=({\bf y})$ which defines the singular locus of $f$ 
is $p+k$, while $q$ is the dimension of the singular set, and $n=p+k+q$.
 The smallest $q$ can be is $p(p+1)/2$. If $q>p(p+1)/2$, then the additional coordinates do not appear in the normal form. For the purposes of this paper, this normal form can be taken as a definition of this type of germ.

\thm 1 Suppose $f:\IC^n\to \IC$ has a non-isolated singularity of type $D(p(p+1)/2,p)$, $n=p(p+1)/2+p$. The homotopy type of the Milnor fiber of $f$ is that of the $S^{2p-1}$ sphere.

\pf We use the technique that appears in the preprint of Fernandez de Bobadilla (\cite{B}). We may chose coordinates so that $f$ is in normal form. Consider the set defined by $f=1$. 
Since $f$ is a homogeneous polynomial its Milnor fiber is diffeomorphic to the set $M$ defined by $f=1$. In turn, $M$ is the total space of the fibration $p$ defined by restricting 
the projection to 
${\bf y}$ to $M$. This is easy to check, as the differential of the map with components $({\bf y},f)$ has maximal rank as long as ${\bf y}\ne 0$, and this holds at all points of $M$. The base of the fibration is 
$(\IC^{p}-0)$, since ${\bf y}\ne 0$. The fiber of $p$ over $({\bf b}) $ in  $\IC^{p}$ consists of the affine hyperplane in $\IC^{p(p+1)/2}$ defined by setting $ {\bf y}$
equal to $({\bf b}) $ in the equation defining $M$. Hence $M$ has a contractible fiber, and has the homotopy type of the base which is the $S^{2p-1}$ sphere.

\cor 2  Suppose $f:\IC^n\to \IC$  has a non-isolated singularity of type $D(q,p)$. The homotopy type of the Milnor fiber of $f$ is that of the $S^{p+n-q-1}$ sphere.

\pf If $n=q+p$, then the Milnor fiber is just the product of the Milnor fiber of the $f$ of Theorem 1.1 with $\IC^{q-p(p+1)/2}$, so the homotopy type doesn't change. 
If $n=q+p+k$, then the normal form for $f$ has $k$ square terms, so the Milnor fiber is the Milnor fiber of the case $n=q+p$ suspended $k$ times, 
which gives a $S^{p+n-q-1}$.

(Fernandez de Bobadilla has advised me by private coorespondence that he also has done this calculation.)

Now we turn to the L\^e numbers.

The L\^e numbers were introduced by Massey, (\cite{Ma}) as a way of relating the Milnor fiber of a function with non-isolated singularities
 to the singularities of the function. For a survey of their properties and more details on their calculation, see \cite{Ma1}. 
They can be calculated in two ways. Let $f:\IC^n\to\IC,0$, $J(f)$ the jacobian ideal of $f$. Consider the blowup
$B_{J(f)}(\IC^n)$ of $\IC^n$ along $J(f)$
 with exceptional divisor $E$. Since $B_{J(f)}(\IC^n)$ lies in $\IC^n\times\IP^{n-1}$, we can intersect both $E$ and 
$B_{J(f)}(\IC^n)$ with $\IC^n\times H$ where $H$ is a linear subspace of $\IP^{n-1}$ We call $H$ a 
plane on $B_{J(f)}(\IC^n).$ If $H_i$ is a generic plane of codimension $i$, then the projection of $B_{J(f)}(\IC^n)\cap H_i$ to $\IC^n$
 is called the relative polar variety of $f$ of codimension $i$, denoted $\Gamma_i(f)$. The projection of $E\cap H_i$ to $\IC^n$
 is a cycle called the L\^e cycle of $f$ of codimension $i+1$. Sometimes we also refer to the cycles which are the components
 of  the L\^e cycle of $f$ of codimension $i+1$ as L\^e cycles as well. The $j$-th cycle of codimension $i$ is then
 denoted $\Lambda_{i+1,j}(f)$. The L\^e number of codimension $i+1$ is the sum of the products of the multiplicity of 
the underlying sets with the degree of the  component cycle.

From this description it is clear that the L\^e cycles can also be constructed by looking at the intersection of $\Gamma_i(f)$ with $S(f)$ and 
calculating the degree of each component of codimension 1 in the polar variety.

The theory of integral closure is tied up with the L\^e cycles, so we briefly discuss it.

Integral dependence is used in local analytic
geometry to relate inequalities between analytic functions, stratification conditions, algebra and analytic invariants. The
basic source for this is the work \cite{LJT} of Monique Lejeune--Jalabert and
 Teissier. (See also \cite {G-2}, \cite{GM}.)

 Let $(X,0)\subset (\IC^N,0)$ be a reduced analytic space germ. Let $I$ be an ideal
in the local ring $\O_{X,0}$ of
$X$ at 0, and $f$ an element in this ring. Then $f$ is integrally dependent
on $I$ if one of the following equivalent conditions obtain:
\itemitem{(i)} There exists a positive integer $k$ and elements $a_j$ in
$I^j,$ so that $f$ satisfies the relation
$f^k+a_1f^{k-1}+\dots+a_{k-1}f+a_k=0$ in
$\O_{X,0}.$
\itemitem{(ii)} There exists a neighborhood $U$ of 0 in $\IC^N,$ a positive
real number $C,$ representatives of the space germ $X,$ the function
germ $f,$ and generators
$g_1,\dots,g_m$ of
$I$ on
$U,$ which we identify with the corresponding germs, so that for all $x$ in
$X$ the following equality obtains:
$$|f(x)|\leq C \max\{|g_1(x)|,\dots,|g_m(x)|\}.$$
\itemitem{(iii)} For all analytic path germs $\phi:(\IC,0)\to(X,0)$ the
pull--back
$\phi^*f$ is contained in the ideal generated by $\phi^*(I)$ in the local ring
of $\IC$ at 0.

\smallskip
If we consider the normalization $\bar B$ of the blowup $B$ of $X$ along the
ideal
$I$ we get another equivalent condition for integral dependence. Denote the
pull--back of the exceptional divisor $D$ of $B$ to $\bar B$ by $\bar D.$
\itemitem{(iv)} For any component $C$ of the underlying set of $\bar D,$ the
order of vanishing of the pullback of $f$ to $\bar B$ along
$C$ is no smaller than the order of the divisor $\bar D$ along $C.$

The elements $f$ in $\O_{X,0}$ that are integrally dependent on $I$ form the
ideal $\bar I,$ the integral closure of $I.$ Often we are only interested in
the properties of the integral closure of an ideal $I;$ so we may replace $I$
by an ideal $J$ contained in $I$ with the same integral closure as $I.$ Such
an ideal $J$ is called a reduction of
$I.$ 

It is easy to see that  $J$ is a reduction of $I$ iff there exists a finite map
$\Bl_IX\to\Bl_JX.$

Now we begin the calculation of the L\^e numbers.

\lem 3  Suppose $f:\IC^n\to \IC$ has a non-isolated singularity of type
$D(p(p+1)/2,p)$, $n=p(p+1)/2+p$. Then $\lambda^{p(p+1)/2-i}(f,0)=0$ for
$i>p>1$.

\pf Let $J$ denote the ideal generated by $(y_1^2,\dots,y_p^2)$. Then $J$ is a reduction of $I^2$. 

To check this just use the curve criterion--given $y_i$ and $y_j$, and a curve $\phi(t)$, then the order of $y_iy_j\circ \phi(t)$, denoted $o(y_iy_j,\phi)$ is greater than or equal to
  min$\{o(y^2_i,\phi),o(y^2_j,\phi)\}$. Since the generators of $I^2$ not in $J$ are in the integral closure of $J$ all of 
$I^2$ is in the integral closure of $J$.

Denote the partial derivatives of $f$ with respect to the $y$  or $x$ variables by $J_y(f)$ or $J_x(f)$. Then from the last paragraph it 
follows that $J_y(f)+J$ is a reduction of $J(f)$.
Now this implies that $B_{J(f)}(\IC^n)$ is finite over $B_{J_y(f)+J}(\IC^n)$; since the second blow-up has fiber over $0$ of dimension at most the number 
of generators of $J_y(f)+J$ less $1$, which is $2p-1$, it follows the same holds for the first blow-up. Now the L\^e cycles are the projection of the intersection of the 
exceptional divisor of $B_{J(f)}(\IC^n)$ with $\IC^n\times H$ where $H$ is a generic plane in $\IP^{n-1}$. So if the codimension of $H$ is $2p$ or more it follows that
$\IC^n\times H$ will miss $E$, hence 
the L\^e cycles of codimension $2p+1$ and more 
must be empty. These are exactly the L\^e cycles of dimension $n-p-i=p(p+1)/2-i$, $i>p>1$.

\cor 4   Suppose $f:\IC^n\to \IC$ has a non-isolated singularity of type
$D(q,p)$. Then $\lambda^{q-i}(f,0)=0$ for $q\ge i>p>1$.

\pf The proof follows the same lines as above noting that $J(f)$ has a reduction with $n+p-q$ elements.

In particular, Corollary 1.4 shows  $\lambda^{0}(f,0)=0$, for $p>1$, so there is no component of the exceptional divisor of $J(f)$ over the origin. In fact, we can describe the components of the exceptional divisor of $J(f)$ completely.

\lem 5 Suppose $f$ has a non-isolated singularity of type $D(q,p)$. Then
$E$ the exceptional divisor of  $B_{J(f)}(\IC^n)$ has only two
components, one of which surjects to
$V(I)$, and the other to $V(I)\cap \det[{\bf X}]$.

\pf  Since   $E$ projects to $V(J(f))=V(I)$, at least one component of $E$ surjects onto $V(I)$, which is smooth, hence irreducible. 
Since on $V(I)$, $f$ has generically an $A(q)$ singularity, only one component surjects.   

At the generic point $V(I)\cap \det[{\bf X}]$, $f$ has the type of a $D(q,1)$ (a Whitney umbrella with 
additional square terms) so again there exists a unique component of $E$ which surjects to
 $V(I)\cap \det[{\bf X}]$. ( That there is exactly one component follows from Proposition 2.2 of \cite {G-3}.)

Now we proceed by induction on $p$, suppose $p>1$ and suppose that there exists $V$, a component of $E$, which maps into $\overline {\Sigma}_2$, 
the sets of points on $V(I)$ where the kernel rank of $[{\bf X}]$ is two or more. Suppose at the generic point of the image of $V$ the kernel rank of  $[{\bf X}]$
is $j$, where $j>1$, since  he kernel rank of $[{\bf X}]$ is two or more. Then since the germs of $f$ along points of constant rank are analytically equivalent, $V$  must surject onto all points of kernel rank $\ge j$.

Picking a suitable transversal $H$ to $\overline {\Sigma}_j$ at the generic point of complementary dimension, the germ of $f|H$ has a singularity of type
$D(j(j+1)/2,j)$, where the dimension of $H$ is $j(j+1)/2+n-q=n-\dim(\Sigma_j)$. Since a component of the exceptional divisor maps to points of kernel rank $j$, the generic fiber dimension over a point of kernel rank $j$ is just the dimension of $V$ less the codimension of the points of kernel rank $j$. Generically, the fiber of the component will be the fiber of  $B_{J(f|H)}(H)$ over the origin, hence the dimension of this fiber is $(n-1)-\dim(\Sigma_j)$, which implies that it is a component of the exceptional divisor of $J(f|H)$. But this contradicts the fact proved in Corollary 1.4, that for $j>1$ there is no component of the exceptional divisor that projects to the origin.

There are two types of L\^e cycles, fixed and moving cycles. The fixed cycles are the images of the components of the exceptional divisor.
 Since the \Af condition holds generically between the
underlying sets of the fixed cycles and the open stratum, a dimension count shows that the 
corresponding component of the exceptional divisor is the conormal of the underlying set. Thus the moving cycles, obtained by intersecting these conormals with 
$\IC^n\times H$ where $H$ is a generic plane in $\IP^{n-1}$ whose codimension is greater than the generic fiber dimension of the conormal over its image,
 are just the polar varieties of the underlying sets of the L\^e cycles.  With this observation we can now calculate the L\^e numbers of $f$. We do this in the simple case of Theorem
1.1, then extend to the general case, which is mostly re-labeling.

\thm 6 Suppose $f:\IC^n\to \IC$ has a non-isolated singularity of type
$D(p(p+1)/2,p)$, $n=p(p+1)/2+p$. Then the L\^e numbers of $f$ are
$$\lambda^{n-p-i}(f)=2^i{{p}\choose p-i } , 0\le i\le p,$$
 $$\lambda^{n-p-i}(f)=0, i>p>1.$$

\pf The second equation is lemma 1.3. The only fixed L\^e cycles are $[V(I)]$ and $2[V(I)\cap \det[{\bf X}]]$. (The coefficients are because a Morse singularity has Milnor number $1$, and a Whitney umbrella has $\lambda^0=2$ and these are the types we get working at a generic point of the underlying set of each cycle, and slicing by a generic transverse plane of complementary dimension.) So $\lambda^{n-p}(f)=1$. 
 Since $V(I)$ is smooth, it has no polar varieties except itself; while the multiplicity at $0$ of $V(I)\cap \det{[\bf X]}$ is $p$, hence
$\lambda^{n-p-1}(f)=2p$.

As usual we assume $f$  in normal form. To compute the remaining L\^e numbers, we first consider the relative polar varieties of $f$. The underlying
 set of the L\^e cycle of dimension $n-p-i$ is gotten by intersecting a relative polar variety of dimension $n-p-i+1$ with $S(f)$.  So consider the set of $p+i-1$
linear combinations of the partial derivatives of $f$, $i>1$. By imposing generic conditions and simplifying our expressions, we can assume they have the form

$$\Pd f{y_l}+q_l(y), p_j(y), 1\le l\le p, 1\le j< i$$

\noindent where $q_l$ and $ p_j$ are quadratic forms in $y$, no $q_i$ in the span of the $p_j$.

   The relative polar variety of $f$ of codimension $p+i-1$ is the closure of the set of points where these expresssions vanish where $y\ne 0$. 
Since all of our expressions are homogeneous, it follows that the underlying sets of the L\^e cycles have a $\IC^*$ action hence are cones. When we projectivize them we get sets of dimension $n-p-i-1$ in 
$\IP^{p(p-1)/2}$. We can compute the multiplicity of the L\^e cycle at the origin by computing the multiplicity of this projective set, ie. its intersection number with a projective space of complementary dimension.

We can slightly change our
our normal form so that we can rewrite the first $p$ equations as
$$2[{\bf X}][{\bf y}] =[{\bf q(y)}].$$                                        

At this point we are going to make a construction which will allow us to find equations for a projective set with the same multiplicity as the projectivization of the underlying set of our L\^e cycle.

Consider the map $G$ from $\IC^n\times \IC\to \IC^n$ given by 
$$G(x,y,\lambda)=(x,\lambda y).$$
Notice that $G$ maps $\IC^n\times 0\to S(f)$. Pullback the $p+i-1$ equations  by $G$.
We get:

$$2[{\bf X}][\lambda{\bf y}] =\lambda^2[{\bf q(y)}]$$
$$\lambda^2{\bf p(y)}=0.$$

The points on the set defined by these equations which map to the relative polar variety are the closure of the subset of points where $\lambda, {\bf y}\ne 0$. These points satisfy the equations:

$$2[{\bf X}][{\bf y}] =\lambda[{\bf q(y)}]$$
$${\bf p(y)}=0.$$

So the points that map to the L\^e cycle are defined by the closure of the points in the above set where $y\ne 0$, and  $\lambda=0$.

This gives the equations:
$$[{\bf X}][{\bf y}] =[0]$$
$${\bf p(y)}=0.$$
The set defined by these equations is a union of components where one component is $\IC^{p(p+1)/2}\times 0$, and the others map to the underlying set of the L\^e cycle. In the components that map to the L\^e cycle, the points where ${\bf y}\ne {\bf 0}$ are dense; further if $({\bf x,y}) $ satisfies the equations and ${\bf y}\ne {\bf 0}$ , then $({\bf x},t{\bf y})$ satisfies the equations for all $t$. Further if $({\bf x})$ is a point on the L\^e cycle, then we may assume that if $({\bf x)}$ is generic, then the kernel rank of the matrix corresponding to $({\bf x})$ is $1$ hence the $({\bf x},t{\bf y})$ is the whole inverse image of such a $({\bf x})$.

These equations define a variety on 
$\IC^{p(p+1)/2}\times \IP^{p-1}$, which maps by projection to $\Lambda^{n-p-i}(f)\subset\IC^{p(p+1)/2}$.

By the observation above, the map is generically 1-1.

Now, consider in $\IP^{p(p+1)/2}\times \IP^{p-1}$ the set of points defined by $[{\bf X}][{\bf y}] =0$, and $\{p_j(y)=0\}$. To calculate the multiplicity at the origin of $\Lambda^{n-p-i}(f)$, it suffices to intersect this last set with $n-p-i-1$ hyperplanes. The number of points realized in $\IP^{p(p+1)/2}\times \IP^{p-1}$ will be the number of lines in the intersection of 
$\Lambda^{n-p-i}(f)$ with the corresponding $n-p-i-1$ hyperplanes in $\IC^n$, and each line counts one toward the multiplicity.

To calculate the number of points, we use the following result from Fulton (\cite{F} p146 eg. 8.4.2)

Given $H_1,\dots, H_{n+m}$ hypersurfaces in $\IP^{n}\times \IP^{m}$ of bidegree $(a_i,b_i)$, then
$$\int [H_1]\dots[H_{n+m}]=\Sigma a_{i_1}\dots a_{i_n}b_{j_1}\dots b_{j_m}$$

\noindent where $i_1<i_2<\dots i_n$ and $j_1<j_2<\dots<j_m$.

In our situation the bi-degree of the equations coming from the matrix are $(1,1)$, while the bidegree of the linear forms is $(1,0)$, and the bidegree of the 
$p_j$ is $(0,2)$. The non-zero terms of the sum are gotten by using for the $j_1,\dots, j_{p-1}$ all $i-1$ terms of bi-degree $(0,2)$, no linear terms and 
$p-i$ terms drawn from the matrix equation. The $i$ part is determined by what remains. The number of such terms is ${{p}\choose {p-i}}$, while the product of the bidegrees is always
$2^{i-1}$. So the multiplicity of the underlying set of $\Lambda^{n-p-i}(f)$ is $2^{i-1}{{p}\choose {p-i}}$. Since the multiplicity of the cycle is $2$, we get
$\lambda^{n-p-i}(f)=2^i{{p}\choose p-i } $, which finishes the proof.

It is interesting to note that the map induced by projection from $\IC^{p(p+1)/2}\times \IP^{p-1}$ to $\IC^{p(p+1)/2}$ when restricted to the incidence variety, $V=\{(x,l)| X|l=0\}$ is a resolution of the variety in $\IC^{p(p+1)/2}$ of symmetric matrices of kernel rank $\ge 1$.

\cor 7 Suppose $f:\IC^n\to \IC$ has a non-isolated singularity of type
$D(q,p)$. Then the L\^e numbers of $f$ are
$$\lambda^{q-i}(f)=2^i{{p}\choose p-i } , 0\le i\le p,$$
 $$\lambda^{q-i}(f)=0, q\ge i>p>1.$$

\pf Let $q_1=q-(p(p+1)/2)$, $k=n-q-p$. Then $q_1$ is the number of coordinates not appearing in the normal form for $f$, and $k$ is the number of coordinates which appear as 
square terms. Assume first $q_1=k=0$. Then we are in the situation of theorem 1.6; rewrite $n-p-i$ as $(p(p+1)/2)-i$. Now keep $k=0$, and increase $q_1$. The effect on the L\^e cycles
is to multiply them by $\IC^{q_1}$, hence the dimensions of the cycles are shifted up by $q_1$, so $n-p-i=(p(p+1)/2)-i$ becomes $(p(p+1)/2)-i+q-1=q-i$. Now let $k$ increase. 
The effect of adding disjoint square terms to the normal form is to leave the L\^e cycles unchanged, so $n-p-i$ can again be replaced by $q-i$.

Massey showed  that the alternating sum of the L\^e numbers of $f$ is the reduced Euler characteristic of the Milnor fiber of $f$. It is a pleasant exercise to recover this result
for the $D(q,p)$ singularities. On the one hand, by Corollary 1.2, the reduced Euler characteristic is $(-1)^{p+n-q-1}(2-1)^p$; expanding $(2-1)^p$, we get
$$(-1)^{p+n-q-1}\Sum_{i=0}^p (-1)^{p-i}2^i{{p}\choose{p-i}}=(-1)^{p+n-q-1}\Sum_{i=0}^p (-1)^{i-p}2^i{{p}\choose{p-i}}$$
$$=\Sum_{i=0}^p (-1)^{(n-1)-(q-i)}\lambda^{q-i}(f). $$

The computations in Theorem 1.6 also compute the polar multiplicities at the origin of the set of $p\times p$ symmetric matrices of kernel rank $\ge 1$.

\cor 8 The polar multiplicities at the zero matrix of the set of $p\times
p$ symmetric matrices of kernel rank $\ge 1$ are given by:

$$m^{p(p+1)/2-i-1}=2^i{{p}\choose{p-i-1}},0\le i<p$$

$$m^{p(p+1)/2-i-1}=0,i\ge p$$

\pf The multiplicities are half the corresponding L\^e numbers, since the L\^e cycles have multiplicity 2.

These computations show another reason why the $D(q,p)$ singularities are interesting. If we put a Whitney stratification on $S(f)$ such that the 
images of the components of the exceptional divisor are a union of strata, then we must include
 many strata which themselves do not correspond to components of the exceptional divisor. These strata are not seen by the exceptional divisor.

Corollary 1.8 shows that the set of $p\times p$ symmetric matrices of kernel rank $\ge 1$ are also interesting--there are many strata in a minimal Whitney stratification which do not correspond
to components of the singular set of the conormal variety of the matrices of kernel rank $\ge 1$. (In fact, a Z-open dense set of the singular set of the conormal variety
 are those points which map to matrices of kernel rank 2.)
 
 This is easy to see as this set is a hypersurface, hence by it is known by \cite{GM} that every component of the conormal modification over the matrices of kernel rank $\ge 2$ has dimension $p(p-3)/2$, yet Corollary 1.8 shows that only over the matrices of kernel rank 2 can the fiber dimension of the conormal be large enough for this to be true.
 
 Now we apply the material of this paper to the calculation of the Euler obstruction of the $D(q,p)$ and of the $p$ by $p$ symmetric matrices.

 The Euler obstruction is an idea introduced by MacPherson (\cite{Mac}) as a
key step in developing the notion of the Chern class for singular spaces. L\^e and Teissier showed that the Euler obstruction of a complex analytic germ at a point $x$ can be computed as the alternating sum of polar multiplicities of $X$ at $x$ (\cite{L-T}). We use this result in our first computation.

\prp 9 The Euler obstruction of the set of $p\times p$ matrices of kernel rank $\ge 1$ is $0$ for $p$ even and $1$ for $p$ odd.

\pf Let $\Sigma_1(p)$ denote  the set of symmetric $p\times p$ matrices of kernel rank $\ge 1$. We know that the reduced Euler characteristic of the Milnor fiber M of a $D(p(p+1)/2,p)$ singularity  satisfies

$$\tilde \chi(M)=(-1)^{2p-1}=-1=\sum_{i=0}^{n-1}(-1)^{n-1-i} \lambda^i(f)$$
$$=(-1)^{p-1}+\sum_{i=n-2p-1}^{n-p-1}(-1)^{n-1-i} 2m^i(\Sigma_1,0).$$

Now there are two cases. If $p$ is even, then we have 
$$0=\sum_{i=n-2p-1}^{n-p-1}(-1)^{n-1-i} m^i(\Sigma_1,0).$$
So the Euler obstruction is $0$.

If $p$ is odd, then we have 
$$-1=\sum_{i=n-2p-1}^{n-p-1}(-1)^{n-1-i} m^i(\Sigma_1,0).$$
$$1=\sum_{i=n-2p-1}^{n-p-1}(-1)^{n-i} m^i(\Sigma_1,0)=Eu(\Sigma_1,0).$$
where the last equality is Cor 5.1.2 \cite {L-T}. (The sum alternates, but the coefficient of the top dimensional term is positive.)

In \cite{B-G} a formula is proved for the Euler obstruction  for $X$ a hypersurface in terms of the L\^e numbers of $f$, where $f$ defines $X$. Using this formula we can easily find the Euler obstruction of a $D(q,p)$ singularity. The formula is:

$${\rm Eu}(X,0)=\chi(L,X,0)+\Sum_{i>0,j}(-1)^{n-i}b(|\Lambda^i_{j,F}(f)|,X){\rm
Eu}(|\Lambda^i_{j,F}(f)|,0).$$

Here $ \Lambda^i_{j,F}(f)$ refers to the $j$-th component of the fixed part of the L\^e cycle of $f$ of dimension $i$, while $b(|\Lambda^i_{j,F}(f)|,X)$ is  the number of spheres in the homotopy type
of the complex link of
$|\Lambda^i_{j,F}(f)|\cap H$ at $z_{i,j}$, where $z_{i,j}$ a generic point of $ \Lambda^i_{j,F}(f)$ and $H$ is a plane of dimension complementary to $|\Lambda^i_{j,F}(f)|$ and transverse to it, while 
$\chi(L,X,0)$ denotes the
Euler
characteristic of the complex link of $X$ at $0$. By a result of Massey's (\cite{Ma2}), this is just
$1+(-1)^{n}m(\Gamma^1(f))$ as 
$m(\Gamma^1(f))$ is the number of spheres in the complex link of $X$ at the origin. Using this formula we can show:

\prp 10 if $f:\IC^n,0\to 0$ has a singularity of type $D(q,p)$ at the origin, $p>1$, $X= f^{-1}(0)$, then 
Eu$(X)=1+(-1)^{n-q}$ if $p$ is even, and $1$ if  $p$ is odd.

\pf We only have 2 fixed L\^e cycles, $S(f)$ and $\Sigma_1$. Further, the relative polar curve is empty for $p>1$ by Cor. 1.7. So we get:

$${\rm Eu}(X)=1+ (-1)^{n-q}+(-1)^{n-q-1}{\rm Eu}(\Sigma_1,0).$$
The second term is the contribution of $S(f)$, while the third term is the contribution of  $\Sigma_1$.
Now the result follows by applying proposition 1.9.

The development of these ideas is continued in \cite{B-G} where the Euler obstruction is computed for those germs which have $A(d)$ or $D(q,p)$ singularities except at the origin. The results of this paper play a key role in applying the formula of \cite{B-G} to these examples.

\sct References

\references
B
J. Fernandez de Bobadilla, {\it Answers to some Equisingularity Questions} preprint 2004

B-G
J. Fernandez de Bobadilla and T. Gaffney, {\it The L\^e numbers of the square of a function and their
applications}, preprint 2005

F
 W. Fulton,
 ``Intersection Theory,''
 Ergebnisse der Mathematik und ihrer Grenzgebiete, 3. Folge
 $\cdot$ Band 2, Springer--Verlag, Berlin, 1984

G-2
   T. Gaffney,
   {\it Integral closure of modules and Whitney equisingularity,}
   \invent 107 1992 301--22
   
   G-3
   T. Gaffney, { \it The multiplicity of pairs of modules and hypersurface
singularities}, Real and Complex Singularities (Sao Carlos, 2004) Trends in Mathematics, Birkhauser 2006, 143-168

G-G
   T. Gaffney and R. Gassler,
      {\it Segre numbers and hypersurface singularities,}
    J. Algebraic Geom. {\bf 8} (1999), 695--736

 GM
   T. Gaffney and D. Massey,
   {\it Trends in equisingularity},  London Math. Soc. Lecture Note Ser.
 {\bf 263} (1999), 207--248

LJT
 M. Lejeune-Jalabert and B. Teissier,
 {\it Cl\^oture integrale des ideaux et equisingularit\'e, chapitre 1}  
Publ.
 Inst. Fourier   (1974)  %\adress St. Martin d'Heres, F--38402

L-T
 D. T. L\^e and B. Teissier, {\it Vari\'et\'es polaires locales et classes de Chern des vari\'et\'es singuli\`eres,}   Ann. of Math. (2)  114  (1981), no. 3,
457--491

Mac
R. D. MacPherson,
{\it Chern classes for singular algebraic varieties,} Ann. of Math. (2)
100 (1974), 423--432

Ma
 D. Massey, The L varieties. I.  Invent. Math.  99  (1990), no. 2, 357--376.

Ma1
  D. Massey {\it L\^e Cycles and Hypersurface Singularities,}
  Springer Lecture Notes in Mathematics 1615 , (1995)

Ma2
D. Massey, {\it Numerical invariants of perverse sheaves}, Duke Math. J.
73 (1994), no. 2, 307--369

P
		 R. Pellikaan,{\it  Hypersurface singularities and resolutions 
of Jacobi modules,} Thesis, Rijkuniversiteit Utrecht, 1985

P1
  R. Pellikaan, {\it Finite determinacy of functions with non-isolated 
singularities,} Proc. London Math. Soc. vol. 57, pp. 1-26, 1988

S1
 D. Siersma, {\it Isolated line singularities,}  Singularities, Part 2 (Arcata, Calif., 1981),  
485--496, Proc. Sympos. Pure Math., 40, Amer. Math. Soc., Providence, RI, 1983

S2
D. Siersma, {\it Singularities with critical locus a $1$-dimensional 
complete intersection and transversal type $A\sb 1$},  Topology Appl.  27  (1987), no. 1, 51--73

\endreferences

\bye